\documentclass{amsart}
\usepackage{amssymb}
\usepackage[dvips]{graphicx}
\usepackage{graphics}
\usepackage{psfrag}
\begin{document}
\newtheorem{thm1}{Theorem}[section]
\newtheorem{lem1}[thm1]{Lemma}
\newtheorem{rem1}[thm1]{Remark}
\newtheorem{def1}[thm1]{Definition}
\newtheorem{cor1}[thm1]{Corollary}
\newtheorem{defn1}[thm1]{Definition}
\newtheorem{prop1}[thm1]{Proposition}
\newtheorem{ex1}[thm1]{Example}
\newtheorem{alg1}[thm1]{Algorithm}


\title[toric ideals of graphs]{Minimal generators of toric ideals of graphs}
\author{Enrique Reyes}
\address{Departamento de Mathem{\' a}ticas,\\ Centro de Investigati{\' o}n y de Estudios Avanzados del IPN\\
Apartado Postal 14-740, 0700 M{\' e}xico City, D.F., M{\' e}xico}
\email{ereyes@math.cinvestav.mx}
\author{Christos Tatakis}
\address{Department of Mathematics, University of Ioannina,
Ioannina 45110, Greece }
\email{chtataki@cc.uoi.gr}
\author{Apostolos Thoma }
\address{Department of Mathematics, University of Ioannina,
Ioannina 45110, Greece }
\email{athoma@uoi.gr}
\thanks{}

%    \subjclass is required.
\subjclass[2000]{Primary 14M25, 05C25, 05C38}
%    The 2010 edition of the Mathematics Subject Classification is
%    now available.  If you are citing a classification from the
%    new scheme, use the following input coding instead.
%\subjclass[2010]{Primary }

\date{}

\dedicatory{}

\begin{abstract}
\par
Let $I_G$ be the toric ideal of a graph $G$. We characterize in
graph theoretical terms primitive, minimal, indispensable and
fundamental binomials of the toric ideal $I_G$.

\end{abstract}
\maketitle


\section{Introduction}

\par  Let $A=\{{\bf a}_1,\ldots,{\bf a}_m\}\subseteq \mathbb{N}^n$
be a vector configuration in $\mathbb{Q}^n$ and
$\mathbb{N}A:=\{l_1{\bf a}_1+\cdots+l_m{\bf a}_m \ | \ l_i \in
\mathbb{N}\}$ the corresponding affine semigroup.  We grade the
polynomial ring $K[x_1,\ldots,x_m]$ over any field $K$ by the
semigroup $\mathbb{N}A$ setting $\deg_{A}(x_i)={\bf a}_i$ for
$i=1,\ldots,m$. For ${\bf u}=(u_1,\ldots,u_m) \in \mathbb{N}^m$,
we define the $A$-{\em degree} of the monomial ${\bf x}^{{\bf
u}}:=x_1^{u_1} \cdots x_m^{u_m}$ to be \[ \deg_{A}({\bf x}^{{\bf
u}}):=u_1{\bf a}_1+\cdots+u_m{\bf a}_m \in \mathbb{N}A.\]  The
{\em toric ideal} $I_{A}$ associated to $A$ is the prime ideal
generated by all the binomials ${\bf x}^{{\bf u}}- {\bf x}^{{\bf
v}}$ such that $\deg_{A}({\bf x}^{{\bf u}})=\deg_{A}({\bf x}^{{\bf
v}})$, see \cite{St}. For such binomials, we define $\deg_A({\bf
x}^{{\bf u}}- {\bf x}^{{\bf v}}):=\deg_{A}({\bf x}^{{\bf u}})$.

Toric ideals have a large number of applications  in several areas
including: algebraic statistics, biology, computer algebra,
computer aided geometric design, dynamical systems,
hypergeometric differential equations, integer programming, mirror
symmetry, toric geometry and graph theory, see \cite{ATY, DS, E-S, MS, St}.
In graph theory there are several monomial or binomial ideals associated to a
graph, see
\cite{CoN, FF, HH, H, NP, SVV1, SVV,  SS, Vi1, Vi, V}, depending on the properties one wishes to study. One of them is the toric
ideal of a graph which has been extensively studied over the last
years, see \cite{CoN, FF, G, G2, K, Ka, OH1, OH, OH2, Hi-O, OH-Ram, VV,
Vi2, Vi}.

The toric ideals are {\em
binomial ideals}, i.e. polynomial ideals generated by binomials. There
are certain binomials in a toric ideal, such as minimal, indispensable, primitive, circuit
and fundamental binomials provide crucial
information about the ideal and therefore they have been studied
in more detail.\\
A binomial $B\in I_A$ is called {\em minimal} if it belongs to at
least one minimal system of generators of $I_A$. The minimal
binomials, up to scalar multiple, are finitely many.
Their number is computed in \cite{ChTh} in terms of combinatorial
invariants of a simplicial complex associated to the toric
ideal. The minimal binomials are characterized as the binomials that can not
be written as a combination of binomials of smaller $A$-degree, see \cite{ChTh,PS}.\\
A binomial $B\in I_A$ is called {\em indispensable} if  there exists a nonzero
constant multiple of it to every
minimal system of generators of $I_A$.
 A recent problem arising from Algebraic Statistics
 is when a toric ideal have a unique minimal
system of binomial generators, see \cite{ChKT, AT}. To study this problem Ohsugi and
Hibi introduced in \cite{Hi-O} the notion of indispensable
binomials and they gave necessary and sufficient conditions for toric
ideals associated with certain finite
graphs to possess unique minimal systems of binomial generators.\\
An irreducible binomial $x^{{\bf u}^+}- x^{{\bf u}^-}$ in $I_{A}$
is called {\em primitive} if there exists no other binomial $
x^{{\bf v}^+}- x^{{\bf v}^-} \in I_{A}$ such that $ x^{{\bf v}^+}$
divides $ x^{{\bf u}^+}$ and $ x^{{\bf v}^-}$ divides $
x^{{\bf u}^-}$.  The set of all primitive binomials  forms the Graver
basis of the toric ideal, see \cite{St}. It follows from the
definition that a non primitive binomial can be written as a sum of products of
monomials times binomials of $I_A$ of smaller
$A$-degree therefore  minimal binomials must be primitive, see also Lemma 3.1 of \cite{OH}.
\\
The support of a monomial $x^{\bf{u}}$ of $K[x_{1},\dots,x_{m}]$
is $supp(x^{\bf{u}}):=\{i\ | \ x_{i}\ divides\ x^{\bf{u}}\}$ and
the support of a binomial $B=x^{\bf{u}}-x^{\bf{v}}$ is
$supp(B):=supp(x^{\bf{u}})\cup supp(x^{\bf{v}})$. An irreducible
binomial $B$ belonging to $I_{A}$ is called a \emph{circuit} of
$I_{A}$ if there is no binomial $B'\in I_{A}$ such that
$supp(B')\subsetneqq supp(B)$. A binomial $B\in I_{A}$ is a
circuit of $I_{A}$ if and only if $I_{A}\cap K[x_{i}\ | i\in
supp(B)]$ is generated by $B$.\\
For a vector ${\bf b}=(b_1,\dots ,b_n)\in \mathbb{N}^n$ we define
$supp({\bf b})=\{i | b_i\neq 0\}$. For a semigroup $\mathbb{N}A$
we denote $K[\mathbb{N}A]$ the semigroup ring of $\mathbb{N}A$.
The semigroup ring $K[\mathbb{N}A]$ is isomorphic to the quotient
$K[x_1,\ldots,x_m]/I_A$, see \cite{MS}. Let ${F}$ be a subset
of $\{1,\dots ,n\}$ then  $A_{{F}}$ is the set
$\{{\bf a}_i| supp({\bf a}_i)\subset {{F}}\}$.  The semigroup
ring $K[\mathbb{N}A_{{F}}]$ is called combinatorial pure
subring of $K[\mathbb{N}A]$, see \cite{OHH} and for a
generalization, see \cite{O}. A binomial $B\in I_A$ is called {\em
fundamental} if there exists a combinatorial pure subring
$K[\mathbb{N}A_{{F}}]$ such that $K[x_i|{\bf a}_i\in A_{{F}}]\cap
I_A=I_{A_{{F}}}=<B>$.

These kinds of binomials are related to each other. The indispensable
binomials are always minimal and the minimal are always primitive.
Also the fundamental binomials are circuits and indispensable, while the
circuits are also primitive. The toric ideals of graphs is the best kind of toric ideals in order to understand
how circuits, fundamentals,
primitive, minimal and indispensable binomials are related, see
Theorems \ref{circuit}, \ref{primitive}, \ref{fundam},
 \ref{minimal}, \ref{indispen}, and to show that the
above relations are strict, see Example \ref{example}.
Actually the toric ideal of a graph gives a way to `view' the
ideal through the graph, but also to construct toric ideals with desired properties.
In the case of the toric ideal of a graph there were several articles in the literature that characterize
these kinds of binomials, most of them for particular cases of graphs, see \cite{G,Ka,OH,OH2,Hi-O,OH-Ram,Vi,V}. 
The aim of this article is to
characterize primitive, minimal, indispensable and fundamental binomials of a toric ideal of a graph
for a general graph and thus understanding better the toric ideal. These characterizations maybe useful
to solve problems in the theory of toric ideals of graphs.

 The results in
this paper are inspired and guided by the work of Oshugi and Hibi
\cite {OH,Hi-O}. In section 2 we present  some
terminology, notations and results about the toric ideals of
graphs. In section 3  we provide the converse of the
characterization of Ohsugi and Hibi \cite{OH} of the primitive
elements of toric ideals of graphs. In section 4 we characterize
the minimal, the indispensable and the fundamental binomials of the toric
ideal of a graph and we give an example that explains the relations
between fundamental, primitive, circuit, minimal and indispensable
binomials. At the end we remark that although the results in the
article are stated and proved for simple graphs, they are also valid with
small adjustments for general graphs with loops and multiple
edges, see Remark \ref{remark}.

\section{Toric ideals of graphs}

\par
In the next chapters,  $G$ will be a finite simple connected graph on the vertex set
$V(G)=\{v_{1},\ldots,v_{n}\}$, except at the final remark \ref{remark} where the graph $G$ may have
multiple edges and loops.  Let $E(G)=\{e_{1},\ldots,e_{m}\}$ be
the set of edges of $G$ and $\mathbb{K}[e_{1},\ldots,e_{m}]$
 the polynomial ring in the $m$ variables $e_{1},\ldots,e_{m}$ over a field $\mathbb{K}$.  We
will associate each edge $e=\{v_{i},v_{j}\}\in E(G)$ with
$a_{e}=v_{i}+v_{j}$ in the free abelian group generated by the
vertices and let $A_{G}=\{a_{e}\ | \ e\in E(G)\}$. With $I_{G}$ we denote
 the toric ideal $I_{A_{G}}$ in
$\mathbb{K}[e_{1},\ldots,e_{m}]$.

\par
A \emph{walk}  connecting $v_{1}\in V(G)$ and
$v_{q+1}\in V(G)$ is a finite sequence of the form
$$w=(\{v_{i_1},v_{i_2}\},\{v_{i_2},v_{i_3}\},\ldots,\{v_{i_q},v_{i_{q+1}}\})$$
with each $e_{i_j}=\{v_{i_j},v_{i_{j+1}}\}\in E(G)$. We call a walk
$w'=(e_{j_{1}},\dots,e_{j_{t}})$ a \emph{subwalk} of $w$ if
$e_{j_1}\cdots e_{j_t}| e_{i_1}\cdots e_{i_q}.$ An edge $e=\{v,u\}$ of a walk $w$
may be denoted also by $(v,u)$ to emphasize the order that the vertices $v, u$ appear
in the walk $w$.
\emph{Length}
of the walk $w$ is called the number $q$ of edges of the walk. An
even (respectively odd) walk is a walk of \emph{even} (respectively odd) length.
A walk
$w=(\{v_{i_1},v_{i_2}\},\{v_{i_2},v_{i_3}\},\ldots,\{v_{i_q},v_{i_{q+1}}\})$
is called \emph{closed} if $v_{i_{q+1}}=v_{i_1}$. A \emph{cycle}
is a closed walk
$$(\{v_{i_1},v_{i_2}\},\{v_{i_2},v_{i_3}\},\ldots,\{v_{i_q},v_{i_{1}}\})$$ with
$v_{i_k}\neq v_{i_j},$ for every $ 1\leq k < j \leq q$. Depending
on the property of the walk that we want to emphasize we may
denote a walk $w$ by a sequence of vertices and edges $(v_{i_1},
e_{i_1}, v_{i_2}, \ldots ,v_{i_q}, e_{i_q}, v_{i_{q+1}})$ or only
vertices $(v_{i_1},v_{i_2},v_{i_3},\ldots,v_{i_{q+1}})$ or only
edges $(e_{i_1},\ldots ,e_{i_q})$  or the edges and vertices that
we want to emphasize and sometimes we separate the walk into
subwalks. For a walk $w=(e_{i_{1}},\dots,e_{i_{s}})$ we denote by
$-w$ the walk $(e_{i_{s}},\dots,e_{i_{1}})$.  Note that, although the graph $G$ has no multiple edges, the
same edge $e$ may appear more than once in a walk. In this case $e$ is
called {\em multiple edge of the walk $w$}. If $w'$ is a subwalk
of $w$ then it follows from the definition of a subwalk that the
multiplicity of an edge in $w'$ is less than or equal to the
multiplicity of the same edge in $w$.

Given an even closed walk $$w =(e_{i_1}, e_{i_2},\cdots,
e_{i_{2q}})$$ of the graph $G$ we denote by
$$E^+(w)=\prod _{k=1}^{q} e_{i_{2k-1}},\ E^-(w)=\prod _{k=1}^{q} e_{i_{2k}}$$
and by $B_w$ the binomial
$$B_w=\prod _{k=1}^{q} e_{i_{2k-1}}-\prod _{k=1}^{q} e_{i_{2k}}$$
belonging to the toric ideal $I_G$. Actually the toric ideal $I_G$
is generated by binomials of this form, see \cite{Vi}. The same
walk can be written in different ways but the corresponding
binomials may differ only in the sign. Note also that different
walks may correspond to the same binomial. For example both walks
$(e_1, e_2, e_3, e_4, e_5, e_6, e_7, e_8, e_9, e_{10})$ and $(e_1,
e_2, e_9, e_8, e_5, e_6, e_7, e_4, e_3, e_{10})$ of the graph $b$
in figure 1 correspond to the same binomial
$e_1e_3e_5e_7e_9-e_2e_4e_6e_8e_{10}$. For convenience by $\bf{w}$
we denote the subgraph of $G$ with vertices the vertices of the
walk and edges the edges of the walk $w$.
If $W$ is a subset of the vertex set $V(G)$ of $G$ then the {\em
induced subgraph} of $G$ on $W$ is the subgraph of $G$ whose
vertex set is $W$ and whose edge set is $\{\{v, u\}\in E(G)|v,u\in
W\}$. When $w$ is a closed walk we denote by $G_w$ the induced
graph of $G$ on the set of vertices  $V({\bf
w})$ of ${\bf w}$. An even
closed walk $w=(e_{i_1}, e_{i_2},\cdots, e_{i_{2q}})$ is said to
be primitive if there exists no even closed subwalk $\xi$ of $w$ of smaller
length such that $E^+(\xi)| E^+(w)$ and $E^-(\xi)| E^-(w)$. The walk $w$
is primitive if and only if the binomial $B_w$ is primitive.

\begin{center}
\psfrag{A}{$e_{1}$}\psfrag{B}{$e_{2}$}\psfrag{C}{$e_{3}$}\psfrag{D}{$e_{4}$}\psfrag{E}{$e_{5}$}
\psfrag{F}{$e_{6}$}\psfrag{G}{$e_{8}$}\psfrag{H}{$e_{7}$}\psfrag{I}{$a$}
\psfrag{K}{$e_{1}$}\psfrag{L}{$e_{2}$}\psfrag{M}{$e_{3}$}\psfrag{N}{$e_{4}$}\psfrag{O}{$e_{5}$}
\psfrag{P}{$e_{6}$}\psfrag{Q}{$e_{7}$}\psfrag{R}{$e_{8}$}\psfrag{S}{$e_{9}$}\psfrag{T}{$e_{10}$}\psfrag{U}{$b$}
\includegraphics{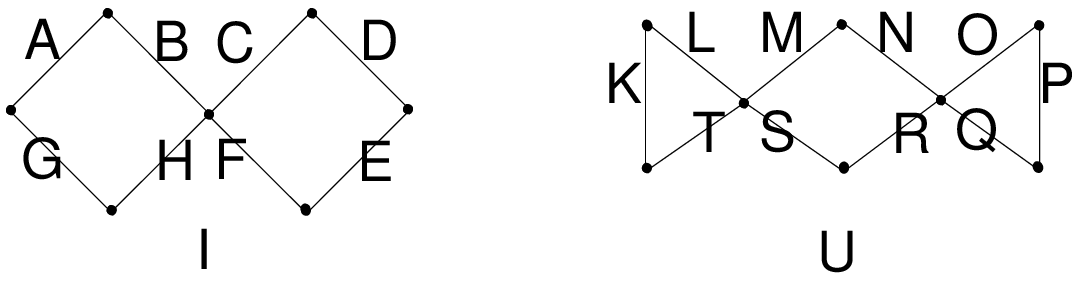}\\
{Figure 1.}
\end{center}

The walk $w=(e_1, e_2, e_3, e_4, e_5, e_6, e_7, e_8)$ of the graph
in Figure $1a$ is not primitive, since there exists a closed even
subwalk of $w$, for example $(e_1, e_2, e_7, e_8)$ such that $e_1
e_7| e_1e_3e_5e_7$ and $e_2e_8|e_2e_4e_6e_8$. While the walk in Figure 1b
$(e_1, e_2, e_3,\\ e_4, e_5, e_6, e_7, e_8, e_9, e_{10})$ is
primitive, although there exists an even closed subwalk $(e_3, e_4,
e_8, e_9)$ but neither $e_3e_8$ divides $e_1e_3e_5e_7e_9$ nor
$e_4e_9$ divides $e_1e_3e_5e_7e_9$.

A necessary characterization of the primitive elements were given by Ohsugi
and Hibi in \cite[Lemma 2.1]{OH}:

\begin{thm1}\label{prim} Let $G$ be a finite connected graph.
If $B\in I_{G}$ is primitive, then we have $B=B_{w}$ where $w$ is
one of the following even closed walks:
\begin{enumerate}
  \item $w$ is an even cycle of $G$
  \item $w=(c_{1},c_{2})$, where $c_{1}$ and $c_{2}$ are odd cycles of $G$ having exactly one common vertex
  \item $w=(c_{1},w_{1},c_{2},w_{2})$, where $c_{1}$ and $c_{2}$ are odd cycles of $G$ having no
  common vertex and where $w_{1}$ and $w_{2}$ are walks of $G$ both of which combine a vertex
  $v_{1}$ of $c_{1}$ and a vertex $v_{2}$ of $c_{2}$.
\end{enumerate}
\end{thm1}
 It is easy to see that any binomial in the first two cases is always
 primitive but this is not true in the third case. Theorem
 \ref{primitive} characterizes completely all primitive binomials.

We will finish this section  with a necessary and sufficient
characterization of circuits that was given by Villarreal in
\cite[Proposition 4.2]{Vi}:
\begin{thm1}\label{circuit}Let $G$ be a finite connected graph. The binomial $B\in I_{G}$ is circuit if and only if $B=B_{w}$ where
\begin{enumerate}
  \item $w$ is an even cycle or
  \item two odd cycles intersecting in exactly one vertex or
  \item two vertex disjoint odd cycles joined by a path.
\end{enumerate}
\end{thm1}

\section{Primitive walks of graphs}

The aim of this chapter is to determine the form of primitive
walks by making more precise the corresponding result  by
Ohsugi-Hibi, see Theorem \ref{prim} or \cite[Lemma 2.1]{OH}.
\\A {\em cut edge} (respectively {\em cut vertex}) is an edge (respectively vertex) of
the graph whose removal increases the number of connected
components of the remaining subgraph.  A graph is called {\em
biconnected} if it is connected and does not contain a cut
vertex. A {\em block} is a maximal biconnected subgraph of a given
graph $G$.
\\Every even primitive walk $w=(e_{i_1},\ldots,e_{i_{2k}})$
partitions the set of edges in the two sets $w^+= \{e_{i_j}|j \
{\it odd}\}, w^-=\{e_{i_j}|j \ {\it even}\}$, otherwise the
binomial $B_w$ is not irreducible.\\
The edges of $w^+$ are called odd edges of the walk and those of
$w^-$ even. Note that for a closed even walk whether an edge is even or
odd depends only on the edge that you start counting from. So it is
not important to identify whether an edge is even or odd but to separate the
edges in the two disjoint classes. \emph{Sink} of a block $B$ is a
common vertex of two odd or two even edges of the walk $w$ which
belong to the block $B$. In particular if $e$ is a cut edge of a
primitive walk then $e$ appears at least twice in the walk and
belongs either to $w^+$ or $w^-$. Therefore both vertices of $e$
are sinks. Sink is a property of the walk $w$ and not of the
underlying graph $\bf{w}$. For example in Figure 1a the walk
$(e_{1},e_{2},e_{7},e_{8})$ has no sink, while in the walk
$(e_{1},e_{2},e_{7},e_{8},e_{1},e_{2},e_{7},e_{8})$ all four
vertices are sinks. Note also that the walk $(e_{1},e_{2},e_{3},e_{4},e_{5},e_{6},e_{7},e_{8})$
in Figure 1a
has one cut
vertex which is not a sink of either block. The walk
$(e_{1},e_{2},e_{3},e_{4},e_{5},e_{6},e_{7},e_{8},e_{9},e_{10})$
in Figure 1b has two cut vertices which are both sinks of all of
their blocks. Theorem \ref{primitive} explains that this is
because the first one is not primitive while the second is.
\begin{center}
\psfrag{D}{$B$}
\includegraphics{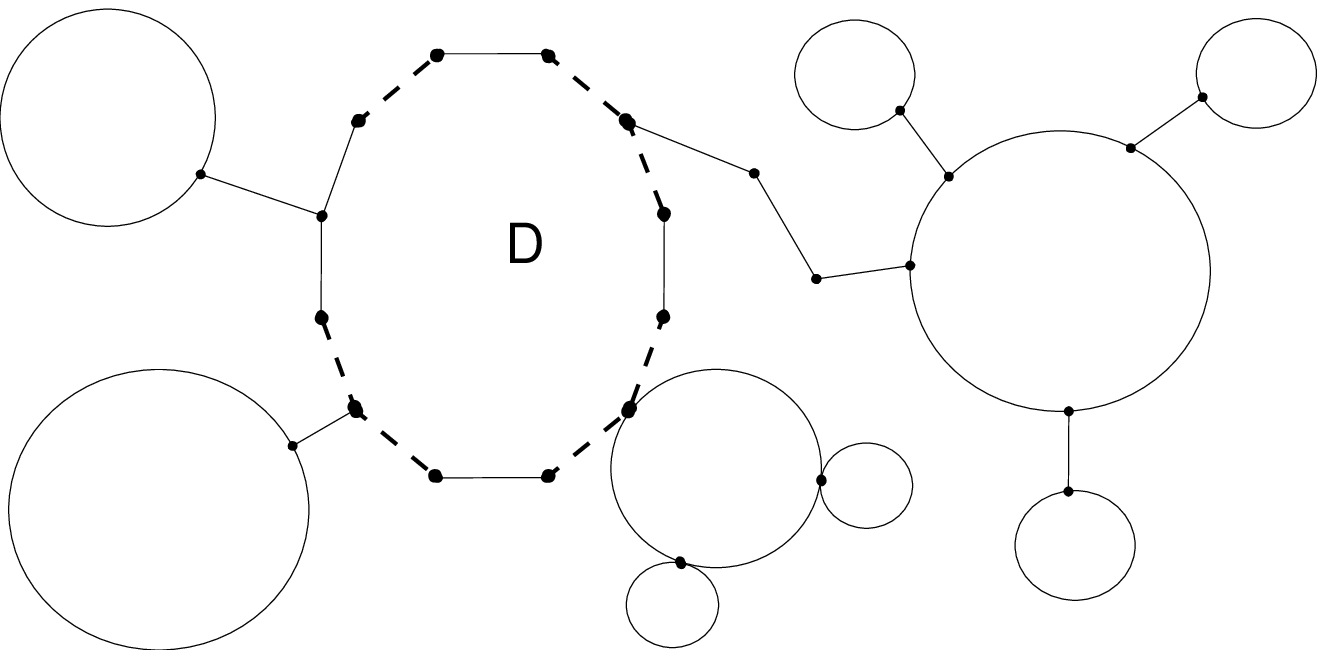}
\\
{Figure 2.}
\end{center}
\begin{thm1} \label{primitive}
Let $G$ a  graph and $w$ an even closed walk of $G$. The walk $w$
is primitive if and only if
\begin{enumerate}
  \item every block of $\bf{w}$ is a cycle or a cut edge,
  \item every multiple edge of the walk $w$ is a double edge of the walk and a cut edge of $\bf{w}$,
  \item every cut vertex of $\bf{w}$ belongs to exactly two blocks and it is a sink of both.
\end{enumerate}

\end{thm1}

\textbf{Proof.} Let $w$ be an even primitive closed walk and let B
be a block of ${\bf w}$ which is not a cut edge. We will prove
that it is a cycle.  Suppose not. Let $w=(e_{i_1},
 \ldots , e_{i_{2s}})$ and $w_B=(e_{i_{j_1}}, \ldots , e_{i_{j_q}})$
 the closed subwalk of
$w$ such that the graph of ${\bf w}_B$ is the block $B$, where
 $e_{i_{j_1}}, \ldots , e_{i_{j_q}}$ are all the edges of the walk $w$ that belong to the the block $B$ and
$j_1<j_2<\dots <j_q$. This is a
closed walk since two blocks intersect in at most one point which is a cut vertex of the graph ${\bf w}$. Since
$B$ is not a cycle, there must be at least one vertex of the walk
$w_B$ which appears twice in $w_B$. If it was exactly one
vertex like that then it should be a cut vertex of $B$
contradicting the biconnectivity of the block $B$. Therefore,  there
must exist at least two vertices $v,u$ of the block $B$ that
appear at least twice in the walk $w_B$ and so $w$ can be written
$w=(v,w_{1},u,w_{2},v,w_{3},u,w_{4})$, where
$w_{1},w_{2},w_{3},w_{4}$ are subwalks of $w$. Note that the
vertices $v,u$ are in this order in $w$ since otherwise $v$ or $u$
will be a cut vertex of $B$. The walk $w$ is primitive, therefore
the closed walk $(v,w_{1},u,w_{2},v)$ is odd, one of the lengths of the subwalks
$w_{1},w_{2}$ has the same parity as the length of $w_{3}$, and both the first edge of
$w_{1}$ and the last of $w_{2}$  belong to $w^+$. Combining
all these, exactly one of the two closed walks
$\xi_1=(v,w_{1},u,-w_{3},v)$ or $\xi_2=(v,w_{3},u,w_{2},v)$ is a
closed even subwalk of $w$ such that either $E^+(\xi_1)| E^+(w)$ and
$E^-(\xi_1)| E^-(w)$ or $E^+(\xi_2)| E^-(w)$ and $E^-(\xi_2)| E^+(w)$.
This contradicts the primitiveness of $w$. So every block is a cycle
or
a cut edge. \\
Let $e=\{u,v\}$ be a multiple edge of $w$. Whenever $e$ appears is either in $w^+$ or $w^-$, since $w$ is a primitive walk.
  The edge $e$ may appear in the walk $w$ in two different
ways, as $(\dots, u,
e, v, \dots)$ or $(\dots, v, e, u, \dots)$. There are two cases.
First case: At least two times the edge appears in the same
way $(\dots, u, e, v, \dots)$ (or $(\dots, v, e, u, \dots)$).
Then the walk $w$ can be written in the form $(u, e, v, w_1, u, e,
v, \dots)$. Since $w$ is primitive and $e$ is written as the first edge of $w$, all the times that $e$ appears
is in $w^+$. Therefore the walk $w_1$ is odd, which means that
$\xi=(u, e, v, w_1, u)$ is an even closed walk, $E^+(\xi)|E^+(w)$
and $E^-(\xi)|E^-(w)$. This contradicts the primitiveness of the walk
$w$. \\
Second case: The edge $e$ appears exactly twice in the walk
and in the two different ways, so $w=(u, e, v, w_1, v, e, u, w_2,
u)$. As before the walks $w_1$, $w_2$ are odd, therefore the first
and the last edges of $w_1$ and $w_2$ all belong to $w^-$.
Suppose that $e$ is not a cut edge of $\bf{w}$ then the $w_1$,
$w_2$ have at least one common vertex $y$. We rewrite $w$ as
$(u,e,v,w_1',y,w_1'',v,e,u,w_2',y,w_2'',u)$. Since $w_2$ is an odd
walk, one of $w_2', w_2''$ is odd and the other is even. Therefore
exactly one of the two walks $(u,e,v,w_1',y,w_2'',u)$,
$(u,e,v,w_1',y,-w_2',u)$ is an even closed walk $\xi$ such that
$E^+(\xi)|E^+(w)$ and $E^-(\xi)|E^-(w)$, contradicting the
primitiveness of the walk $w$. We conclude that $e$ is a double
edge of the walk $w$ and a cut edge of $\bf{w}$.\\
Let $v$ be a cut vertex, then it belongs to at least two blocks.
Since $v$ is a cut vertex  $w$ can be written as $w=(v, e_1,\dots
,e_s, v, e_{s+1}, \dots, e_t, v, \dots )$. where $e_{1}, e_{s}$
are in the same block $B$ and $\{e_i|1\leq i\leq s\}\cap \{e_i|s+1\leq i\leq t\}=\emptyset$. Then $e_{1},e_{s}$ are both in
$w^+$. Otherwise $(v, e_{1},\dots,e_{s},v)$ is an even
closed subwalk of $w$, contradicting the primitiveness of the walk
$w$. So $v$ is a sink and  the subwalk $(v,e_{1},\dots,e_{s},v)$  is odd. Similarly
the walk $(v,e_{s+1},\dots,e_{t},v)$ is odd and $e_{s+1}, e_t$ are  both in $w^-$.
Then $w'=(v, e_1,\dots ,e_s, v, e_{s+1}, \dots, e_t, v)$ is
an even subwalk of $w$ such that $E^+(w')|E^+(w)$ and $E^-(w')|E^-(w)$
and since $w$ is primitive $w'=w$. We
conclude that $v$ belongs to exactly two blocks of ${\bf w}$ and
it is a sink of both.
\newline
Conversely let $w$ be an even closed walk satisfying the three
conditions of the Theorem which is not primitive. Then there exists
a primitive subwalk $w'$ of $w$ of smaller length than $w$, such
that $E^+(w')|E^+(w)$ and $E^-(w')|E^-(w)$. From the first part of
the proof we know that also $w'$ satisfies  the three conditions
 of the Theorem \ref{primitive}. We claim that the graphs ${\bf w}$ and ${\bf w}'$
have exactly the same blocks.
Let $B_{w'}$ be a block of $w'$ then there exists a block $B_w$ of
$w$ such that $B_{w'}\subset B_w$. From the first condition
$B_{w'}$ is a cut edge or a cycle. Suppose that $B_{w'}=\{e\}$ is
a cut edge of
${\bf w}'$ then $e$ must be double edge of $w'$.
Since $E^+(w')|E^+(w)$ and $E^-(w')|E^-(w)$ the edge
$e$ is a multiple edge of $w$ and therefore from the second
condition a cut edge of ${\bf w}$, thus a block of ${\bf w}$. In
the  case that $B_{w'}$ is a cycle obviously $B_w$ is the same cycle and
therefore $B_{w'}= B_w$. So all blocks of $\mathbf{w}'$ are blocks
of $\mathbf{w}$. Conversely suppose that there exist a block of $\mathbf{w}$
which is not a block of $\mathbf{w}'$. Since ${\bf w}$ is connected
 there must be at least
one block of $\mathbf{w}$ which is not a block of $\mathbf{w}'$
and has a contact point with $\mathbf{w}'$. Then this point should
be a sink of both since $E(w')^+| E(w)^+$ and $E(w')^-| E(w)^-$.
But if it is a sink of $w'$ then it should belong to exactly two
blocks of ${\bf w}'$. This implies that it should belong to at least three
blocks of ${\bf w}$, a contradiction to the third property of $w$.\\ Therefore the graphs ${\bf w}$ and ${\bf
w}'$ are identical and every simple edge of the walk $w'$ is a
simple edge of $w$ and every double edge (cut edge) of the walk
$w'$  is
 a double edge of $w$. Therefore $E^+(w')=E^+(w)$ and
 $E^-(w')=E^-(w)$. Therefore they have the same length,
a contradiction. We conclude that $w$ is primitive. \hfill $\square$

From Theorem \ref{primitive} easily follows the following corrolary that describes the underlying graph
of a primitive walk.

\begin{cor1} \label{primitive-graph}
Let $G$ a  graph and $W$ a subgraph of $G$. The subgraph $W$ is the graph  ${\bf w}$ of a primitive walk $w$
 if and only if
\begin{enumerate}
  \item every block of $W$ is a cycle or a cut edge and
  \item every cut vertex of $W$ belongs to exactly two blocks and separates the graph in two parts, the total number of edges
of the cyclic blocks in each part is odd.
\end{enumerate}

\end{cor1}

\section{Minimal and indispensable binomials of graphs  }

The first aim of this section is to characterize the walks $w$ of the graph $G$ such that the binomial $B_w$ belongs to a minimal system
of generators of the ideal $I_G$. Certainly the walk has to be primitive, but this is not enouph. The walk must have  more properties,
the first one it depends on the graph ${\bf w}$ and the rest on the induced graph $G_w$ of $w$, see Proposition \ref{sprim} and Theorem \ref{minimal}.
\begin{def1}
A binomial $B\in I_{G}$ is called minimal if it belongs to a
minimal system of binomial generators of $I_{G}$.
\end{def1}

\begin{def1} We call strongly primitive walk a primitive walk that has
not two sinks with distance one in any cyclic block.
\end{def1}

\begin{prop1}\label{sprim} Let $w$ be an even closed walk such that the binomial $B_w$ is minimal
then the walk $w$ is strongly primitive.
\end{prop1}
\textbf{Proof.}

\begin{center}
\psfrag{A}{$\xi_1$}\psfrag{B}{$\xi_3$}\psfrag{D}{$\xi_2$}\psfrag{C}{$e$}\psfrag{E}{$u$}\psfrag{F}{$v$}
\includegraphics{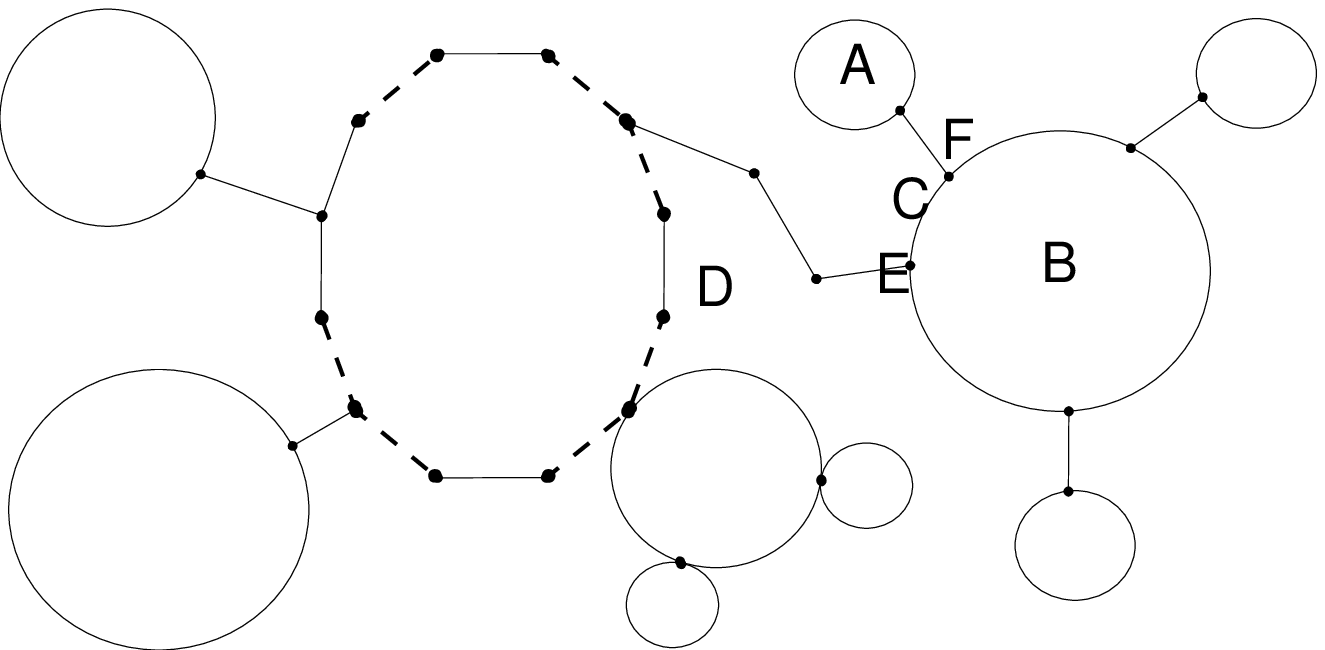}\\
{Figure 3.}
\end{center}

 The binomial $B_{w}$ is minimal therefore the walk $w$ is primitive.
Suppose that $w$ is not strongly primitive, then there exist two
sinks $v, u$ of the same block $B$ with distance one. We will call
$e$ the edge $\{v,u\}$. Then $w$ can be written as
$(v,\xi_{1},v,e,u,\xi_{2},u,\xi_{3},v)$ for some walks $\xi_1, \xi_2, \xi_3$, where at least the first and the
last edge of $\xi_{3}$ are in the block $B$. The walks $\xi_1,\ \xi_2$ are
closed walks and since $w$ is primitive they are necessarily odd.
Therefore $\xi_3$ is also odd. So the first and the last edge of
$\xi_{3}$, as also $e$ are in $w^-$. Since $v,u$ are sinks the closed
walks $w_{1}=(\xi_{1},e,\xi_{2},e)$ and $w_{2}=(e,\xi_{3})$ are
both even and the binomial
$B_{w}=B_{w_{1}}\frac{E^+(w_{2})}{e}+B_{w_{2}}\frac{E^-(w_{1})}{e}$
is not minimal, a contradiction. Note that $E^+(w_2)/e\neq 1\neq
E^-(w_1)/e$, otherwise the even closed walk $w_1$ or $w_2$ has length 2,
which is impossible since  $G$ has no multiple edges. \hfill $\square$

While the property of a walk to be primitive depends only on the
graph ${\bf w}$, the property of the walk to be minimal or
indispensable depends also on the induced graph $G_w$. An edge $f$
of the graph $G$ is called a \emph{chord} of the walk $w$ if the vertices of the edge
$f$ belong to $V(\bf{w})$ and $f\notin E (\bf{w})$. In other words an
edge is called chord of the walk $w$ if it belongs to $E(G_w)$ but not in $E(\textbf{w})$. Let $w$ be an
even closed walk
$((v_{1},v_{2}),(v_{2},v_{3}),\ldots,(v_{2k},v_{1}))$ and
$f=\{v_{i},v_{j}\}$ a chord of $w$. Then $f$ \emph{breaks} $w$ in two
walks:
$$w_{1}=(e_{1},\ldots,e_{i-1}, f, e_{j},\ldots,e_{2k})$$ and
$$w_{2}=(e_{i},\ldots,e_{j-1},f),$$ where $e_{s}=(v_{s},v_{s+1}),\ 1\leq s\leq 2k$ and $e_{2k}=(v_{2k},v_{1}).$
The two walks are both even or both odd.

In the next definition we are interested in chords of the walk. We
partition the set of chords of a primitive even walk in three
parts: bridges, even chords and odd chords.
\begin{def1} A chord $f=\{v_{1},v_{2}\}$
is called bridge of a primitive walk $w$  if
there exist two different blocks $B_{1},B_{2}$ of $\bf{w}$ such that
$v_{1}\in B_{1}$ and $v_{2}\in B_{2}$. A chord is called even
(respectively odd) if it is not a bridge and breaks the walk in
two even walks (respectively in two odd walks).
\end{def1}

In the walk of Figure 4, there are three chords which are bridges
of $w$, those marked by $b$ and there is one chord which is even,
it is marked by $c$. In the walks of Figure 5, all chords
are odd. Note that the two vertices of a bridge may also belong to the same block,
for example that happens in one of the three bridges in Figure 4.

\begin{center}
\psfrag{A}{$i$}\psfrag{B}{$j$}\psfrag{C}{$c$}\psfrag{D}{$B_{s}$}\psfrag{E}{$b$}\psfrag{F}{$b$}\psfrag{G}{$b$}
\includegraphics{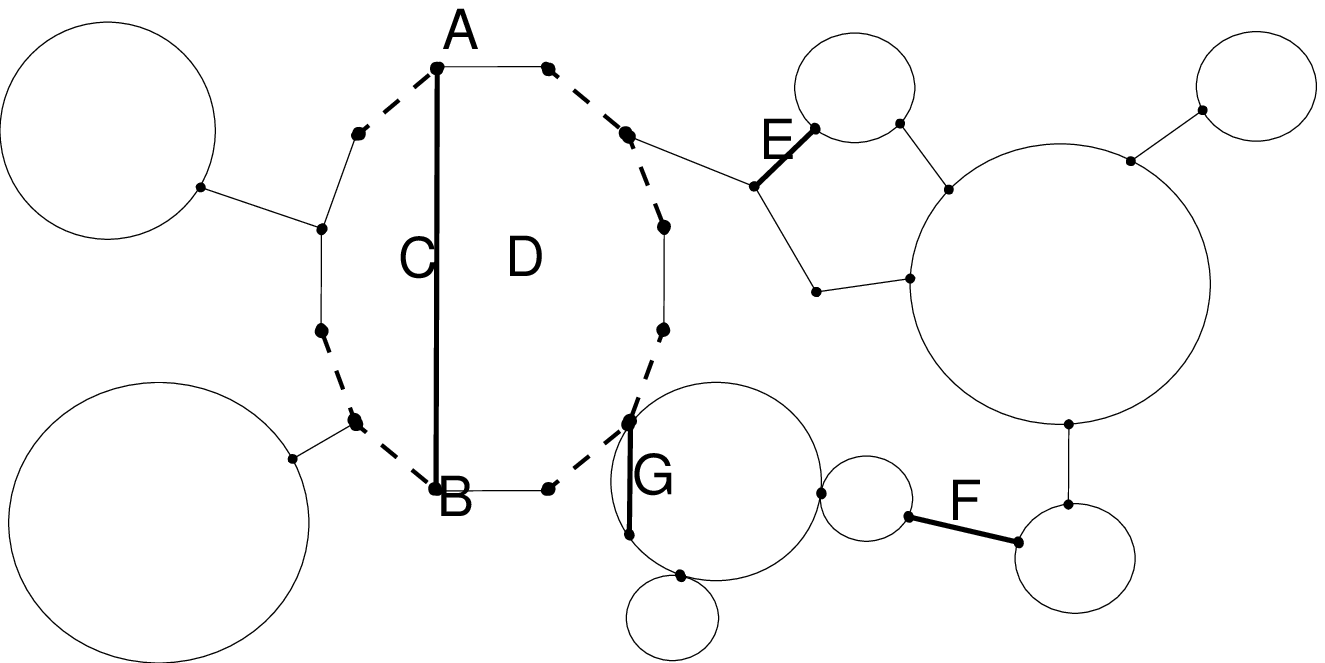}\\
{Figure 4.}
\end{center}

The next definition generalizes the corresponding definitions of
Ohsugi and Hibi, see \cite{Hi-O}.
\begin{def1}\label{creffec} Let
$w=((v_{i_{1}},v_{i_{2}}), (v_{i_{2}},v_{i_{3}}),\cdots ,
(v_{i_{2q}},v_{i_{1}}))$ be a primitive walk. Let
$f=\{v_{i_{s}},v_{i_{j}}\}$ and $f'=\{v_{i_{s'}},v_{i_{j'}}\}$ be two
odd chords (that means not bridges and $j-s,j'-s'$ are even) with $1\leq
s<j\leq 2q$ and $1\leq s'<j'\leq 2q$. We say that $f$ and $f'$
cross effectively in $w$ if $s'-s$ is odd (then necessarily $j-s',
j'-j, j'-s$ are odd) and  either $s<s'<j<j'$ or $s'<s<j'<j$.
\end{def1}

\begin{def1}\label{F_4}
We call an $F_4$ of the walk $w$ a cycle $(e, f, e', f')$ of
length four which consists of two edges $e,e'$ of the walk $w$
both odd or both even, and two odd chords  $f$ and $f'$ which
cross  effectively in $w$.
\end{def1}
In Figure 5 there are two cyclic blocks of primitive walks and in
each one exactly two odd chords which cross effectively. In the
first block they form an $F_{4}$, while in the second they do not.
Combining Definitions \ref{creffec} and \ref{F_4} two odd chords
are part of an $F_4$ if $i'-j=\pm 1$ and $j'-i=\pm 1$, or
$i'-i=\pm 1$ and $j'-j=\pm 1$.
\begin{center}
\psfrag{A}{$f$}\psfrag{B}{$f'$}\psfrag{C}{$B_{s}$}\psfrag{D}{$F_{4}$}
\psfrag{E}{$Cross \ effectively\ odd \ chords$}\psfrag{F}{$f'$}\psfrag{G}{$f$}
\includegraphics{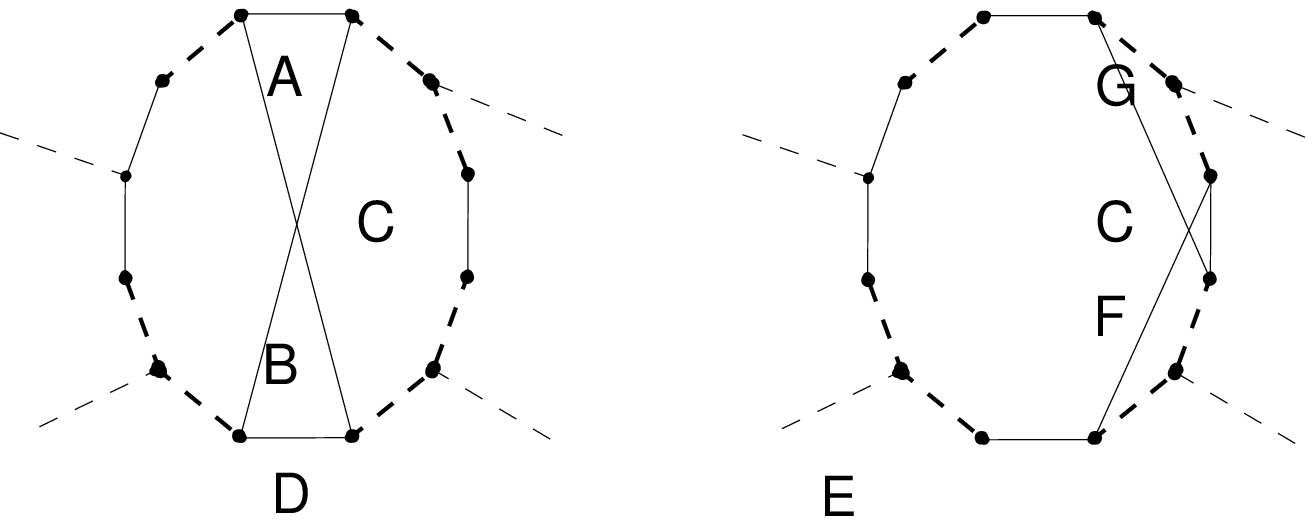}\\
{Figure 5.}
\end{center}

\begin{def1}\label{stronglycross}
Let $w$ be a primitive walk and $f, f'$ be two odd chords. We say that $f, f'$ cross strongly effectively
in $w$ if they cross effectively  and they do not form an $F_4$ in $w$.
\end{def1}

\begin{prop1} \label{minimal1} Let $w$ be a primitive walk. If $B_{w}$ is a minimal binomial then
 all the chords of $w$ are
odd and there are not two of them which cross  strongly effectively.
\end{prop1}
\textbf{Proof.} Let $w=(e_1, e_2, \dots ,e_{2s})$ be a
primitive walk. If $B_{w}$ is a minimal binomial, then from
Proposition \ref{sprim} it follows that $w$ is strongly primitive.
Let $e=\{v_1, v_{2l-1}\}$ be an even chord of $w$, and let
$w_1=(e_1, e_2, \dots ,e_{2l-1},e)$, $w_2=(e, e_{2l},\dots
,e_{2s})$ be the two even walks that $e$ breaks $w$. Then
$B_{w}=B_{w_{1}}\frac{E^+(w_{2})}{e}-B_{w_{2}}\frac{E^-(w_{1})}{e}$,
so $B_w$ is not minimal. Note that $E^+(w_2)/e\neq 1\neq
E^-(w_1)/e$, since $G$ has no multiple edges. \\Suppose that a
minimal binomial $B_w$ has a bridge $e=\{v_1, v_2\}$. Since $v_1,
v_2$ belong to different blocks there must be at least one cut
vertex $v$ such that the walk $w$ can be written $(v, w_1, v_1,
w_2, v, w_3, v_2, w_4, v)$. Note that if $v=v_{1}$ or $v=v_{2}$
one of the walks $w_{1}$, $w_{4}$ is empty. The closed walks $(v, w_1,
v_1, w_2, v)$ and $(v, w_3, v_2, w_4, v)$ are both odd, otherwise
$B_w$ is not primitive. Therefore one of the $w_1, w_2$ has to be
odd and the other even. Similarly for $w_3, w_4$. Note also that the four
walks $(v,w_{1},v_{1},w_{2},v,w_{3},v_{2},w_{4},v)$ ,
$(v,w_{1},v_{1},w_{2},v,-w_{4},v_{2},-w_{3},v)$,
$(v,-w_{2},v_{1},-w_{1},v,w_{3},v_{2},w_{4},v)$ and
$(v,-w_{2},v_{1},-w_{1},v,-w_{4},v_{2},-w_{3},v)$ give the same
binomial. Therefore we can assume that $w_{1},w_{3}$ are odd and
$w_{2},w_{4}$ are even. Then the two closed walks
$\zeta_{1}=(w_{2},w_{3},e)$ and $\zeta_{2}=(w_{4},w_{1},e)$ are
even and
$B_{w}=B_{\zeta_{1}}\frac{E^+(\zeta_{2})}{e}-B_{\zeta_{2}}\frac{E^-(\zeta_{1})}{e}$
is not minimal, a contradiction. Note that
$\frac{E^+(\zeta_{2})}{e}\neq 1 \neq
\frac{E^-(\zeta_{1})}{e}$, since $G$ has no multiple edges.\\
Suppose now that $w$ has two odd chords  $f=\{v_1, v_2\} ,f'=\{u_1, u_2\}$ which cross strongly
effectively in
$w$. Then $w$ is in the form
$(v_1, w_1, u_1, w_2, v_2, w_3, u_2, w_4, v_1)$.
 We have that the walks $\xi_{1}=(w_{1},f',-w_{3},f)$, and
$\xi_{2}=(w_{2},f,-w_{4},f')$ are even, since the walks
$(w_{1},w_{2},f),(w_{2},w_{3},f'),(w_{3},w_{4},f), (w_{4},w_{1},f')$ are
odd. Then
$B_{w}=B_{\xi_{1}}\frac{E^+(\xi_{2})}{ff'}-B_{\xi_{2}}\frac{E^-(\xi_{1})}{ff'}$
is not minimal, a contradiction. Note that since the odd chords $f, f'$ do not form an
$F_{4}$, $\frac{E^+(\xi_{2})}{ff'}\neq 1 \neq \frac{E^-(\xi_{1})}{ff'}$.
\hfill $\square$

\begin{def1} Two primitive walks $w, w'$ differ by an
$F_4$, $\xi=(e_1,f_1,e_2,f_2)$, if $w=(w_1, e_1, w_2, e_2)$ and $w'=(w_1,
f_1,-w_2, f_2)$, where both $w_1, w_2$ are odd walks. Two even closed walks $w, w'$
are $F_4$-equivalent if either $w=w'$ or there exists a series of walks $w_1=w, w_2, \dots, w_{n-1}, w_n=w'$ such that
$w_i$ and $w_{i+1}$ differ by an $F_4$, where $1\leq i\leq n-1$.
\end{def1}

Note that if $w$ and $w'$ are $F_4$-equivalent then the induced graphs $G_w$ and $G_{w'}$ are
equal. We denote by $L_w$ the equivalence class of $w$ under the $F_4$-equivalent relation.

\begin{prop1} \label{minimal2} If the primitive walks $w$ and $w'$
are $F_4$-equivalent then $B_w$ is minimal if and only if $B_{w'}$ is minimal.
\end{prop1}
\textbf{Proof.} Suppose that $w=(w_1, e_1, w_2, e_2)$ and $w'=(w_1,
f_1,-w_2, f_2)$ are even closed walks which differ by an $F_{4}$, where $F_{4}$ is $\xi=(e_1,f_1,e_2,f_2)$.
Then $B_{w}=B_{w'}-\frac{E^-(w)}{e_{1}e_{2}}B_{\xi}$ and the result follows.\hfill $\square$

The $F_4$ separates the vertices of ${\bf w}$ in two parts $V({\bf w}_1), V({\bf w}_2)$, since both edges $e_1, e_2$  of an $F_4$, $(e_1, f_1, e_2, f_2)$, belong to the same block of $w=(w_1, e_1, w_2, e_2)$,

\begin{def1} We say that an odd chord $f$ of a primitive walk $w=(w_1, e_1, w_2, e_2)$ crosses an $F_4$, $(e_1, f_1, e_2, f_2)$,
if  one of the vertices of $f$ is in $V({\bf w}_1)$, the other in $V({\bf w}_2)$ and $f$ is different from  $f_1, f_2$.

\end{def1}

\begin{prop1} \label{minimal3} Let $w$ be a primitive walk. If $B_{w}$ is a minimal binomial, then
 no odd chord crosses an $F_4$ of the walk $w$.
\end{prop1}
\textbf{Proof.} Let $B_w$ be a minimal binomial. Suppose that there exists an odd chord $f=\{v_1, v_2\}$ that crosses the $F_4$, $(e_1, f_1, e_2, f_2)$, of the walk
$w=(w_1, e_1, w_2, e_2)$.  Then $w$ can be written in the form $(w_1',v_1,w_1'', e_1, w_2',v_2, w_2'', e_2)$. The chord $f$ is odd
therefore the walks $(f, w_2'',e_2, w_1')$ and $(f, w_1'',e_1, w_2')$ are both odd. Also, since $(e_1, f_1, e_2, f_2)$ is an $F_4$,
the walks $w_1$ and $w_2$ are both odd. Therefore $(w_1'', f_1, -w_2'',f)$ and $(w_1', f, -w_2',f_2)$ are both even.
So, from the definition, $f$ is an even chord of $w'=(w_1, f_1, w_2, f_2)$. Note that $f$ is not a bridge of $w'$ since it is not a bridge
of $w$. Therefore from Proposition \ref{minimal1} $B_{w'}$ is not minimal and from Proposition \ref{minimal2} $B_w$ is not minimal, a contradiction.
\hfill $\square$

In fact, for a primitive walk $w$  the  walks in $L_w$ are primitive, an $F_4$ of $w$ is an $F_4$ for all walks in $L_w$, although sometimes chords and edges change role.
A bridge of $w$ is a bridge for every walk in $L_w$ and odd chords of $w$ (respectively even chords)
are odd chords (respectively even chords) for every walk in $L_w$, except if they cross an $F_4$. In the last case they may change parity, it depends on how many $F_4$
they cross.

\begin{thm1} \label{minimal}
Let $w$ be an even closed walk. $B_{w}$ is a minimal binomial if
and only if $w$ is strongly primitive, all the chords of $w$ are
odd, there are not two of them which cross  strongly effectively and no odd chord crosses an $F_4$ of the walk $w$.
\end{thm1}
\textbf{Proof.} The one direction follows from Propositions \ref{minimal1} and \ref{minimal3}.  \newline For the converse, let $w$ be an
even closed walk such that all the chords of $w$ are odd, there
are not two of them which cross strongly effectively and no odd chord crosses an $F_4$ of the walk $w$.
Suppose that $B_{w}$ is not minimal. Then there exists a minimal
walk $\delta$ such that $E^+(\delta) | E^+(w)$ and
$E^+(\delta) \not= E^+(w)$,
thus edges of $\delta^{+}$ are edges of $w^{+}$.
 We have
$\deg_{A_G}(E^-(\delta))=\deg_{A_G}(E^+(\delta))<\deg_{A_G}(E^+(w))=\deg_{A_G}(E^-(w))$. This means that the vertices of $\delta^{-}$ are in ${\bf w}$
and so edges of
$\delta^{-}$ are edges or chords of $w$, which means actually they are odd
chords by hypothesis.\\ We claim that every such $\delta $ is an $F_4$ of $w$. Suppose not, then among all those
walks $\delta $ which are not $F_4$ of $w$ and  $E^+(\delta) | E^+(w)$ and
$E^+(\delta) \neq E^+(w)$, we choose one, $\gamma$, such that $\gamma $ has the fewest possible chords of $w$.\\
First case: The walk $\gamma $ does not have any chords of $w$, then all edges of $\gamma^{-}$ are edges of $w$, so
$\gamma^{+}\subset w^{+}$ and $\gamma^{-}\subset w$ and since $w$
is primitive then there exists at least one $e \in \gamma^{-}\cap
w^{+}$. Therefore $\gamma=(\dots, e_{1},e,e_{2},\dots)$, where all edges
$e_{1},e,e_{2}$ are  in $w^{+}$. Note that whenever
there are two blocks joined by a cut vertex, the adjoining edges in the two different blocks have different parity, since the walk $w$ is primitive.
Thus all the edges $e_{1},e,e_{2}$  are in one block
of $w$, which necessarily is a cycle and then the two vertices in
between are sinks of $w$. A contradiction to strongly
primitiveness. Note that if two of $e_{1},e,e_{2}$ are the same edge, then
this edge will be a double edge of $w$ and therefore a cut edge of ${\bf w}$, so the edges $e_{1},e,e_{2}$ are in two blocks, a contradiction.\\
Second case: $\gamma^{-}$ has at least one chord of $w$. Then
$\gamma=(w_{1},f_{1},w_{2},f_{2},\dots,w_{s},f_{s})$ where
$w_{1},\dots,w_{s}$ are subwalks of $w$ and $f_{1},\dots,f_{s}$
are odd chords of $w$ satisfying the hypotheses and $s$ is minimal. Both vertices of
an odd chord $f$ of $w$ are in the same cyclic block, thus $f$ divides
$w$ into two regions $w^{+}(f),w^{-}(f)$. There must exist at
least one chord $f_{i}$  such that the region $w^{+}(f_{i})$ does
not contain a chord.
The last edge of $w_i$ and the first of $w_{i+1}$ are in $\gamma^+\subset w^+$. The chord $f_i$ is odd, therefore the
one of these two edges is in $w^+(f_i)$ and the other in $w^-(f_i)$. Without loss of generality we can suppose that the first edge of $w_{i+1}$
is in $w^+(f_i)$. The walk $\gamma$ is closed and none of the vertices of $f_i$ is a cut vertex of $\gamma $, since $f_i$ is not a bridge of $w$, therefore
 there must be a chord which has a
vertex in $w^{+}(f_{i})$ and a vertex in $w^{-}(f_{i})$. This chord is   the $f_{i+1}$ since $w^+(f_i)$ does not contain a chord.
Let $f_i=\{v_{i_{s}},v_{i_{j}}\}$ and $f_{i+1}=\{v_{i_{s'}},v_{i_{j'}}\}$.
Since the
first and the last edge of $w_{i+1}$ are in $\gamma^{+}\subset
w^{+}$, $s'-j$ (the number of edges of $w_{i+1}$) is odd. But from
the hypothesis $f_{i},f_{i+1}$ can not cross effectively except if
they form an $F_{4}$, which means that either $|s'-j|=|j'-s|=1$ or
$|s'-s|=|j'-j|=1$.  In the first case $w$ is an $F_4$, the $(e_{i_s}, f_i, e_{i_{s'}}, f_{i+1})$.
In the second case there exists an even minimal walk
$\gamma'=(w_{1},f_{1}, \dots ,w_{i},e_{i_{s+1}},-w_{i+1},e_{i_{s+1}},w_{i+2},\dots)$
with two less chords, a contradiction to the minimality of the chords  of $\gamma$.

We conclude that if for a walk $\delta $ we have  $E^+(\delta) | E^+(w)$ and
$E^+(\delta) \neq E^+(w)$, then $\delta$ is an $F_4$ of $w$. Remark that the conditions of Theorem \ref{minimal} if they are satisfied by
the walk $w$, then they are also satisfied by any other walk in $L_w$. We fix a minimal set $\{B_{w_1},\cdots ,B_{w_t}\}$ of binomial generators for the ideal $I_{A_G}$,
 for some even closed walks $w_1, \dots ,w_t$. For a
 $w' \in { L}_w$ we define
$$
r(w')=min\{\sum_{l=1}^t\vert g_l \vert  \mid  B_{w'}=\sum_{i=1}^k
g_l B_{w_l}  \}
$$
and $\vert g_l \vert$ is the number of monomials of $g_l$. We take a walk $v \in
{ L}_w$ such that $r(v)$ is minimal. We claim that $B_v$ is one of the minimal generators $B_{w_1},\cdots ,B_{w_t}$.
Suppose not,
then it is written in the form
$B_{v}=E^+(v)-E^-(v)=\sum_{r=1}^q g_{i_r}B_{w_{i_r}}$ and without loss of generality we can suppose that
 $E^+(w_{i_1})|E^+(v)$ and
$E^+(v)/E^+(w_{i_1})\not=1 $ is a monomial in $g_1$. Then  $w_{i_1} $ is necessarily an $F_4$, $(e_1, f_1, e_2, f_2)$, of $v=(v_1, e_1, v_2, e_2)$ and $e_1e_2=E^+(w_{i_1})|E^+(v)$.
Consider $v'=(v_1, f_1, v_2, f_2)$, then $v'\in L_w$ and $$B_{v'}=E^+(v')-E^-(v')=\frac{E^+(v)}{e_{1}e_{2}}f_1f_2-E^-(v)=(g_1-\frac{E^+(v)}{e_{1}e_{2}})B_{w_{i_1}}+\sum_{r=2}^q
g_{i_r}B_{w_{i_r}}. $$ But then  $r(v')<r(v)$ which is a contradiction. Note that the coefficients of the monomials in $g_i$
are $1$ or $-1$, see \cite{ChKT, DS, St} for more information about the generation of a toric ideal.\\ Therefore $B_v$ is minimal and
from Proposition \ref{minimal2} $B_w$ is
minimal. \hfill $\square$

Note that in the  cases in  Theorem \ref{minimal} where we have
more than one $F_4$,  two $F_4$ of the walk cannot have a common edge
 and they cannot cross, since in all these cases we get
an odd chord which  crosses an $F_4$.

\begin{thm1}\label{indispen}
Let $w$ be an even closed walk. $B_{w}$ is an
indispensable binomial if and only if $w$ is a strongly primitive walk, all the chords of $w$ are
odd and there are not two of them which cross effectively.
\end{thm1}
\textbf{Proof.} From Proposition \ref{minimal2} if $w$ has an $F_4$ then $B_w$ is not indispensable, since it can
be replaced by $B_{w'}$. So $w$ has not an $F_4$, then the result follows from Theorem \ref{minimal}.
 \newline
Conversely, let $w$ be a strongly primitive walk, all the chords of $w$ are
odd and there are not two of them which cross effectively.
Suppose that $B_{w}$ is not indispensable. Then there exists a minimal
walk $\delta \neq w$ such that $E^+(\delta) | E^+(w)$,
thus edges of $\delta^{+}$ are edges of $w^{+}$.
 We have
$\deg_{A_G}(E^-(\delta))=\deg_{A_G}(E^+(\delta))\leq \deg_{A_G}(E^+(w))=\deg_{A_G}(E^-(w))$. This means that the vertices of $\delta^{-}$ are in ${\bf w}$
and so edges of
$\delta^{-}$ are edges or chords of $w$, which means actually they are odd
chords by hypothesis.\\
First case: The walk $\delta $ does not have any chords of $w$. In that case the proof is exactly the same as in the corresponding
part in the proof of Theorem \ref{minimal}.\\
Second case: $\delta^{-}$ has at least one chord of $w$. Then
$\delta=(w_{1},f_{1},w_{2},f_{2},\dots,w_{s},f_{s})$ where
$w_{1},\dots,w_{s}$ are subwalks of $w$ and $f_{1},\dots,f_{s}$
are odd chords of $w$ satisfying the hypotheses. There must exist at
least one chord $f_{i}$  such that the region $w^{+}(f_{i})$ does
not contain a chord.
The last edge of $w_i$ and the first of $w_{i+1}$ are in $\delta^+\subset w^+$. The chord $f_i$ is odd, therefore the
one of these two edges is in $w^+(f_i)$ and the other in $w^-(f_i)$. Without loss of generality
we can suppose that the first edge of $w_{i+1}$
is in $w^+(f_i)$. The walk $\delta$ is closed and none of the vertices of $f_i$ is a cut vertex of $\delta $, since $f_i$ is not a bridge of $w$, so
 there must be a chord of $w$ which has a
vertex in $w^{+}(f_{i})$ and a vertex in $w^{-}(f_{i})$. This chord is   the $f_{i+1}$ since $w^+(f_i)$ does not contain a chord.
Let $f_i=\{v_{i_{s}},v_{i_{j}}\}$ and $f_{i+1}=\{v_{i_{s'}},v_{i_{j'}}\}$.
Since the
first and the last edge of $w_{i+1}$ are in $\delta^{+}\subset
w^{+}$,  the number of edges $(s'-j)$ of $w_{i+1}$ is odd. Therefore $f_i$ and $f_{i+1}$ cross effectively, a contradiction.
\\ Therefore $B_w$ is indispensable. \hfill $\square$

Remark that combining Theorem \ref{minimal}, Proposition \ref{minimal2} and Theorem \ref{indispen},
we have that if $B_{w}$ is indispensable then $w$ has no $F_{4}$ and if $B_{w}$ is minimal but not
indispensable then $B_{w}$ has at least one $F_{4}$. If no minimal generator has an $F_4$ then the toric ideal
is generated by indispensable binomials, so the ideal $I_G$ has a unique system of
binomial generators and conversely.

An even closed walk $w$ of a graph $G$ is called {\em
fundamental} if for every even closed walk $w'$ of the induced
subgraph of $G_{w}$ it holds $B_{w'}\in <B_w>$. A binomial
$B_w$ is fundamental if $w$ is fundamental, see \cite{Hi-O}.

\begin{thm1}\label{fundam} If $w$ is an even closed walk, then the binomial
$B_w$ is fundamental if and only if $w$ is a circuit and
has no chords except in the case that it is a cycle with no even chords and at most one odd chord.
\end{thm1}
\textbf{Proof.} Let $w$ be an even closed walk such that the binomial
$B_{w}$ is fundamental. From \cite[Theorem 1.1.]{OH2} we know that
$B_{w}$ is a circuit and $B_{w}$ is an indispensable binomial.
Since $B_{w}$ is a circuit, from Proposition \ref{circuit} there
are three cases. If $w$ is a cycle the result follows from
\cite[Lemma 4.2.]{OH2}. In the other two cases, $w$ is a circuit
with no even chords and bridges, since $B_{w}$ is indispensable.
Suppose that  $w$ has an odd chord. The odd chord necessarily is a
chord of  one of the two odd cycles. Every chord of an odd cycle
breaks the cycle in two cycles, one of which is odd and the other even. The
even cycle gives a binomial in $I_{G_w}$ which is not in $<B_w>$.
A contradiction arises since $B_{w}$ is fundamental.
\newline Conversely if $w$ is a cycle with no even chords and at
most one odd chord, the result follows from \cite[Lemma
4.2.]{OH2}. On the other hand if $w$ is not a cycle then it is a
circuit with no chords. Therefore $w$ has no even cycles and
$B_{w}$ is fundamental. \hfill $\square$

\begin{ex1}\label{example}{\rm The simplest possible graph which
shows that the relations between fundamental, primitive,
indispensable, minimal binomials and circuits are strict is the
following: let $G$ be the graph with 10 vertices and 14 edges of
figure 6.

\begin{center}
\psfrag{A}{$e_{1}$}\psfrag{B}{$e_{2}$}\psfrag{C}{$e_{3}$}\psfrag{D}{$e_{4}$}\psfrag{E}{$e_{5}$}\psfrag{F}{$e_{6}$}\psfrag{G}{$e_{7}$}
\psfrag{H}{$e_{8}$}\psfrag{I}{$e_{9}$}\psfrag{J}{$e_{10}$}\psfrag{K}{$e_{11}$}\psfrag{L}{$e_{12}$}\psfrag{M}{$e_{13}$}\psfrag{N}{$e_{14}$}
\includegraphics{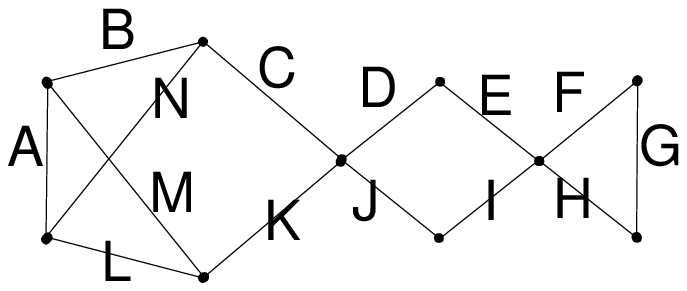}\\
{Figure 6.}
\end{center}

The Graver basis has twenty two elements:
$B_1=e_2e_{12}-e_{13}e_{14}$, $B_2=e_2e_{11}-e_3e_{13}$,
$B_3=e_3e_{12}-e_{11}e_{14}$, $B_4=e_4e_9-e_5e_{10}$,
$B_5=e_{14}e_{2}e_{4}^2e_{6}e_{8}-e_{1}e_{3}^2e_{5}^2e_{7}$,
$B_6=e_{14}e_{2}e_{10}^2e_{6}e_{8}-e_{1}e_{3}^2e_{9}^2e_{7}$,
$B_7=e_{13}e_{12}e_{4}^2e_{6}e_{8}-e_{1}e_{11}^2e_{5}^2e_{7}$,
$B_8=e_{13}e_{12}e_{10}^2e_{6}e_{8}-e_{1}e_{11}^2e_{9}^2e_{7}$,
$B_9=e_{14}e_{2}e_{4}e_{6}e_{8}e_{10}-e_{1}e_{3}^2e_{5}e_{7}e_{9}$,
$B_{10}=e_{13}e_{12}e_{4}e_{6}e_{8}e_{10}-e_{1}e_{11}^2e_{5}e_{7}e_{9}$,
$B_{11}=e_{3}e_{1}e_{11}e_{5}^2e_{7}-e_{2}e_{12}e_{4}^2e_{6}e_{8}$,
$B_{12}=e_{3}e_{1}e_{11}e_{9}^2e_{7}-e_{2}e_{12}e_{10}^2e_{6}e_{8}$,
$B_{13}=e_{3}e_{1}e_{11}e_{5}e_{9}e_{7}-e_{2}e_{12}e_{4}e_{10}e_{6}e_{8}$,
$B_{14}=e_{3}e_{1}e_{11}e_{5}^2e_{7}-e_{14}e_{13}e_{4}^2e_{6}e_{8}$,
$B_{15}=e_{3}e_{1}e_{11}e_{9}^2e_{7}-e_{14}e_{13}e_{10}^2e_{6}e_{8}$,
$B_{16}=e_{3}e_{1}e_{11}e_{5}e_{9}e_{7}-e_{14}e_{13}e_{4}e_{10}e_{6}e_{8}$,
$B_{17}=e_{14}e_{1}e_{11}^2e_{9}^2e_{7}-e_{2}e_{12}^2e_{10}^2e_{6}e_{8}$,
$B_{18}=e_{14}e_{1}e_{11}^2e_{5}^2e_{7}-e_{2}e_{12}^2e_{4}^2e_{6}e_{8}$,
$B_{19}=e_{14}e_{1}e_{11}^2e_{9}e_{5}e_{7}-e_{2}e_{12}^2e_{4}e_{10}e_{6}e_{8}$,
$B_{20}=e_{13}e_{1}e_{3}^2e_{9}^2e_{7}-e_{12}e_{2}^2e_{10}^2e_{6}e_{8}$,
$B_{21}=e_{13}e_{1}e_{3}^2e_{5}^2e_{7}-e_{12}e_{2}^2e_{4}^2e_{6}e_{8}$,
$B_{22}=e_{13}e_{1}e_{3}^2e_{9}e_{5}e_{7}-e_{12}e_{2}^2e_{4}e_{10}e_{6}e_{8}$.
The first eight of them are fundamental binomials. The first ten
are indispensable binomials and the first sixteen binomials are
minimal. Note that the number of minimal generators $\mu(I_G)$ is
$13$ and there are $8$ different, up to non zero constants,
minimal systems of binomial generators. The cause of the
dispensability of the binomials $ B_{11},\ldots,B_{16}$ is the
existence of an $F_4$, $(e_2,e_{13},e_{12}, e_{14})$. The cause of
the primitive elements $B_{17},B_{18},B_{19}$ and
$B_{20},B_{21},B_{22}$ not to be minimal is the existence of
bridges: $e_3$ in the first three and $e_{11}$ in the last three. Finally
all of them are circuits except the binomials
$B_{9},B_{10},B_{13},B_{16},B_{19},B_{22}$. }
\end{ex1}

\begin{rem1}\label{remark}{\rm
For simplicity of the statements and the proofs we assumed that
the graphs are simple. But actually most of the results are valid
with small adjustments for graphs with loops and multiple edges.
Theorem \ref{primitive} about primitive walks is true exactly as
it is stated, but note that you may have cycles with one edge, a
loop, and cycles with two edges, in the case that you have
multiple edges between two vertices. The property of a walk to
give a minimal binomial depends on the induced graph and  one may
have chords which are loops or multiple edges. In this case in
Theorem \ref{minimal}, which describes the even closed walks that
determine minimal generators, chords which are multiple edges are
also permitted and loops where the vertex of the loop is not a cut
vertex of ${\bf w}$. While in Theorem \ref{indispen} chords which
are multiple edges are not permitted, but chords which are loops such that the
vertices of the loops are not cut vertices of ${\bf w}$ are permitted.

}
\end{rem1}

\end{document}